\documentclass[preprint,12pt]{elsarticle}

\usepackage{amsthm}
\usepackage{endnotes}
\let\footnote=\endnote

\usepackage{algorithmic}
\usepackage{algorithm}
\usepackage{amsmath,amsfonts,amssymb}
\usepackage{caption}
\usepackage{subcaption}
\usepackage{soul}
\usepackage[official]{eurosym}
\usepackage{xcolor}
\usepackage{hyperref}

\newtheorem{thm}{Theorem}
\newtheorem{lemma}{Lemma}
\newtheorem{corollary}{Corollary}

\newtheorem{prop}[thm]{Proposition}
\newtheorem{definition}{Definition}

\begin{document}
\begin{frontmatter}




\title{An exact dynamic programming approach to segmented isotonic regression}


\author[1]{V\'ictor Bucarey}
\ead{vbucarey@vub.be}
\author[2,5]{Martine Labb\'e}
\ead{mlabbe@ulb.ac.be}
\author[3]{Juan M. Morales\corref{cor1}}
\ead{juan.morales@uma.es}
\author[4]{Salvador Pineda}
\ead{spinedamorente@gmail.com}

\cortext[cor1]{Corresponding author}
\address[1]{Data Analytics Laboratory, Vrije Universiteit Brussel, Brussels, Belgium}
\address[2]{D\'epartement d'Informatique, Universit\'e Libre de Bruxelles, Brussels, Belgium}
\address[3]{Department of Applied Mathematics, University of Malaga, Malaga, Spain}
\address[4]{Department of Electrical Engineering, University of
Malaga, Malaga, Spain}
\address[5]{INRIA Lille Nord-Europe, France}


\begin{abstract}
This paper proposes a polynomial-time algorithm to construct the monotone stepwise curve that minimizes the sum of squared errors with respect to a given cloud of data points. The fitted curve is also constrained on the maximum number of steps it can be composed of and on the minimum step length. Our algorithm relies on dynamic programming and is built on the basis that said curve-fitting task can be tackled as a shortest-path type of problem. 
Numerical results on synthetic and realistic data sets reveal that our algorithm is able to provide the globally optimal monotone stepwise curve fit for samples with thousands of data points in less than a few hours. Furthermore, the algorithm gives a certificate on the optimality gap of any incumbent solution it generates. 
From a practical standpoint, this piece of research is motivated by the roll-out of smart grids and the increasing role played by the small flexible consumption of electricity in the large-scale integration of renewable energy sources into current power systems. Within this context, our algorithm constitutes an useful tool to generate \emph{bidding} curves for a pool of small flexible consumers to partake in wholesale electricity markets. 
\end{abstract}


\end{frontmatter}


%


\section{Introduction}
In this paper, we deal with the problem of how to fit a curve to a given cloud of data points under the conditions that the fitted curve must be non-increasing (or non-decreasing) and piecewise constant (or, equivalently, \emph{stepwise}), with a predefined limited number of pieces (also referred to as \emph{steps} or \emph{blocks} in what follows). This problem is inspired by the bidding rules that large consumers or a pool of small consumers must comply with when participating in an electricity market. Their bids for purchasing electricity in these markets must be often submitted in the form of a non-increasing stepwise price-consumption curve, for which the maximum number of bid blocks is also constrained. These curves reflect how consumers value electricity and therefore, their sensitivity to its price (which is referred to as consumers' elasticity), see, for instance, \citet{su2009quantifying}. Furthermore, beyond its use for market bidding, the consumers' sensitivity to the electricity price constitutes essential information for the design of tariff schemes and demand response programs (\citet{grimm2020optimal,soares2020designing,zugno2013bilevel}). Indeed, with the advent of Information and Communications Technologies and the roll-out of the so-called \emph{smart grids}, small consumers of electricity are being provided with the means to actively adjust their consumption in response to the electricity price. However, their consumption patterns are still uncertain, dynamic and affected by other factors different from the electricity price. The result is that estimating a bidding curve that properly reflects consumers' price-sensitivity is a statistical challenge. This paper provides an algorithm to efficiently compute that curve from a set of price-consumption observations.

Beyond the practical context that inspires this piece of research, our work is closely related to various thrusts of research or thematic areas that also motivate it, namely:
\begin{description}
    \item[\textbf{Statistical regression}.] We desire to fit a monotonically decreasing curve  to a given cloud of data points, while satisfying the following two extra conditions: i) The fitted curve must be piecewise constant and ii) there is a maximum number of pieces the fitted curve can be comprised of.
    While the literature review includes a wealth of research papers analyzing related concepts and tools such as \emph{isotonic} regression (see, e.g., \citet{mair2009isotone}, \citet{Tibshirani2011}, 
    and references therein), \emph{segmented} regression (\citet{muggeo2003estimating}), and the popular \emph{multivariate adaptive regression spline} (\citet{friedman1991multivariate}), these regression techniques produce fitted curves that fail to satisfy at least one of the conditions mentioned above. Furthermore, they are frequently based on iterative, greedy or heuristic algorithms. 
    Indeed, the fitted response of the isotonic regression is a monotone piecewise constant function (although efficient algorithms to produce smooth continuous functions are also available, see, e.g., \citet{sysoev2019smoothed}), but is not limited in the number of pieces it may be comprised of. For its part, segmented regression leads to curve fits that are not necessarily monotone.
    Against this background, we propose an exact shortest-path algorithm that is capable of delivering, in polynomial time, the monotone stepwise curve (with a maximum of $K$ steps) that constitutes the globally optimal data fit according to the least-squares criterion.
    
    We remark that, as pinpointed in \citet{lerman1980fitting}, the stepwise shape of the target curve releases the fitting process from the continuity condition at the breakpoints that is typically enforced in segmented regression, thus making it computationally easier. On the other hand, we additionally impose that the fitted curve be non-increasing, which adds an extra layer of complexity to the regression problem at hand. Actually, to our knowledge, the works that are the closest to ours are those of \citet{hawkins1976point} and \citet{dahl2000cardinality}. In the former, they describe a dynamic programming approach to perform segmented regression over a sequence of observations with at most $K$ segments and no continuity requirement at the transition points. \citet{dahl2000cardinality} offer an interesting computational perspective on this same problem, which they pose as a cardinality-constrained shortest path problem over a restricted class of acyclic-directed graphs known as \emph{2-graphs} and for which they propose several solution algorithms.  Our approach, in contrast, works with generic acyclic-directed graphs (which do not need to be 2-graphs), that is, we allow for arcs between any pair of nodes $i$ and $j$, with the only condition that $i>j$ in the topological order induced by the acyclicity of the graph. Furthermore, neither \citet{hawkins1976point}, nor \citet{dahl2000cardinality} consider any monotonicity constraint, which is, though, critical to our problem (seen as an extension or generalization of isotonic regression) and to the practical application that motivates it. 
    
    Finally, we mention that \citet{rote2018isotonic} uses dynamic programming for isotonic regression, but, again, with no constraint on the number of pieces  the fitted curve can be made up of.
    
    \item[\textbf{Inverse optimization}.] Recently, inverse optimization has emerged as a promising mathematical framework to infer the input parameters to an optimization problem that have given rise to a series of optimal (or quasi-optimal) solutions (\citet{ahuja2001inverse,esfahani2018data, chan2019inverse}). In the last few years, inverse optimization has been widely used to infer consumers' utility from a certain product (\citet{keshavarz2011imputing,aswani2018inverse}), in particular, electricity (\citet{saez2016data, saez2017short}). Essentially, it is often assumed that the market behavior of a pool of (rational) electricity consumers is driven by the following maximization problem 
    \begin{equation*}
    \underset{x \geq 0}{\text{Maximize}} \enskip \int_{0}^{x}{b(s)ds} - p x
\end{equation*}
where $p x$ is the payment the pool of consumers has to make for purchasing $x$ units of electricity in the market at price $p$, and $b(\cdot)$ is the so-called \emph{bidding curve} expressing the response of the consumers to the electricity price. Many electricity markets around the world require that this bidding curve be non-increasing and stepwise, with a maximum number of steps. Dealing with this problem by way of inverse optimization involves estimating the step function values of this curve and the breakpoints from a series of observed pairs $\{(\hat{p}_i,\hat{x}_i)\}_{i=1}^{I}$. As highlighted in \citet{aswani2018inverse}, however, the available estimation approaches based on inverse optimization may result in statistically \emph{inconsistent} estimators or require the reformulation of the problem as a bilevel (NP-hard) problem. In this regard, \citet{aswani2018inverse} propose a statistically consistent polynomial-time semiparametric algorithm to tackle a certain class of inverse optimization problems. Nevertheless, the regression problem we address here, when seen from the lens of inverse optimization, does not comply with the conditions that ensure the statistical and polynomial-time performance of their algorithm, because some of the parameters to be estimated, specifically, the breakpoints, appear in the constraints defining the feasible region of the forward problem. In contrast, we propose an algorithm that directly solves the statistically consistent formulation of the problem to optimality in polynomial time. 

\item[\textbf{Unsupervised learning}.] The problem we address in this paper can be also interpreted as a clustering problem through which a series of observed pairs $\{(\hat{p}_i,\hat{x}_i)\}_{i=1}^{I}$ are grouped in such a way that: 
\begin{enumerate}
\item There is a maximum number of clusters $K$ into which the data points can be grouped into.
\item The resulting clusters must satisfy some connectivity constraints. In our particular case, these connectivity constraints impose that only clusters with adjacent prices $\hat{p}_i$ can be merged together (see, e.g., \citet{guo2009greedy}). 
\item If $(p^*_m,x^*_m)$ and $(p^*_n,x^*_n)$ are the centroids of clusters $m$ and $n$, respectively, then $x^*_m \geq x^*_n \iff p^*_m \leq  p^*_n$ in order to guarantee a non-increasing curve. 
\end{enumerate}
The technical literature includes some works in which structured clustering is used in power system applications. For instance, 
\citet{pineda2018chronological} propose a hierarchical clustering methodology to approximate  time series that are used to determine the optimal expansion planning of the European electricity network. Due to the usual NP-hard nature of clustering methods, the clusters are often obtained through computationally efficient greedy algorithms. However, to the best of our knowledge, the technical literature does not report any clustering methodology that simultaneously satisfies the three conditions specified above. Therefore, our work also contributes to the realm of structured data clustering. 
\end{description}

The rest of this paper is organized as follows. In Section~\ref{PF}, we formulate the curve-fitting problem that we aim to solve. Section~\ref{s:shortest_algorithm} introduces the solution algorithm we propose to that end, which is based on dynamic programming and, more specifically, on the cardinality-constrained shortest path problem. Section~\ref{AS} provides various strategies to accelerate said algorithm, whose performance is subsequently tested in Section~\ref{NE} using synthetic data sets and a data set coming from  a real-life practical application. Lastly, conclusions are duly drawn in Section~\ref{Conc}. 

\section{Problem definition}\label{PF}
Consider a given set of pairs of points on the real plane $\{(\hat{p}_i,\hat{x}_i)\}_{i=1}^{I}$. Without loss of generality, we assume that $\hat{p}_1 < \hat{p}_2 < \ldots < \hat{p}_I$, while the set of indexed coordinates $\{\hat{x}_i\}_{i=1}^{I}$ may not exhibit any particular order. Let $\mathcal{F}$ be the class of real functions $f:[\hat{p}_1,\hat{p}_I] \rightarrow\mathbb{R}$ 
that are non-increasing and piecewise constants, with at most $K$  blocks or steps, $K\in \mathbb{Z}_+$. We seek to solve the following least-square minimization problem, hereinafter referred to as \emph{LSP}:
\begin{equation}
\text{(LSP)}  \hspace{2cm}  \min_{f \in \mathcal{F}}  \enskip \sum_{i=1}^{I}{\left(\hat{x}_i-f(\hat{p}_i)\right)^2} \hspace{2cm}\label{eq:ERM} 
\end{equation}
A function $f$ member of the class $\mathcal{F}$ can be expressed as 

\begin{equation}
f(p) = \sum_{k=1}^{K}{u_k \mathbb{I}_{[p_k,p_{k+1})}(p)} 
\end{equation}
where $\mathbb{I}_{[p_k,p_{k+1})}(p)$ is the indicator function equal to 1 if $p_k\leq p < p_{k+1}$, and 0 otherwise. Again without loss of generality, we set $p_1=\hat{p}_1$ and use ${p}_{K+1} = \hat{p}_{I+1} > \hat{p}_{I}$ as a dummy $p$-coordinate to guarantee that all $\hat{p}_i$ are covered by the solution. Besides, $u_1\geqslant u_2 \geqslant \ldots \geqslant u_K \geqslant 0$  represent the \emph{step values} of the blocks. We remark that functions  $f \in \mathcal{F}$ with less than $K$ blocks can be also represented in this way, since two consecutive blocks are allowed to have the same function value. 

Using this characterization of the class of functions $\mathcal{F}$ and taking $p_1=\hat{p}_1$  and $p_{K+1} > \hat{p}_{I}$ a dummy price coordinate as mentioned above, problem~\eqref{eq:ERM} can be recast as follows.
\begin{subequations}\label{eq:ERM_2}
\begin{align}
    \min\limits_{\textbf{u},\textbf{p}} & \enskip \sum_{i=1}^{I}{\left(\hat{x}_i-\sum_{k=1}^{K}{u_k \mathbb{I}_{[p_k,p_{k+1})}   (\hat{p}_i)}\right)^2}\label{eq:ERM_2_OF}\\
    \textrm{s.t.}& \enskip u_k \geqslant u_{k+1}, \enskip \forall k \leqslant K-1\label{eq:ERM_2_monotonicity}\\
    \phantom{\textrm{s.t.}}& \enskip p_{k+1} \geqslant p_{k}, \enskip \forall k \leqslant K 
\end{align}
\end{subequations}

Determining the breakpoints $\{p_k\}_{k=2}^K$, which are needed to compute the indicator functions appearing in the objective function~\eqref{eq:ERM_2_OF}, constitutes the major source of complexity in problem~\eqref{eq:ERM_2}. Constraint~\eqref{eq:ERM_2_monotonicity}, which enforces the non-increasing character of the fitted curve, also adds another layer of difficulty to the selection of those breakpoints. The easiest task in problem~\eqref{eq:ERM_2} is to compute the values $u_k$ that minimize the squared error~\eqref{eq:ERM_2_OF} for a given set of  intervals $[p_k,p_{k+1})$. 
In the following section, we introduce a shortest path algorithm through which we can solve problem~\eqref{eq:ERM_2} in polynomial time.  

\section{Resource constrained shortest-path algorithm}
\label{s:shortest_algorithm}

We begin by demonstrating that problem~\eqref{eq:ERM_2} (and hence, problem LSP) can be equivalently reformulated as a cardinality-constrained shortest-path problem. The equivalence stems from the evidence that the optimal breakpoints are within the $p$-coordinates of the cloud of points $\{(\hat{p}_i,\hat{x}_i)\}_{i=1}^{I}$.

\subsection{Properties and problem reformulation}
\label{ss:properties}

For a function $f \in \mathcal{F}$, the objective function value of Problem (\ref{eq:ERM_2}) can be rewritten as:
\begin{align}
 \sum_{k=1}^K \sum_{i: p_k \leq \hat{p}_i < p_{k+1}} (\hat{x}_i - u_k)^2\label{eq:OBJ}
 \end{align}
 
The following lemma shows that we can restrict the search of breakpoints to the set of $p$-coordinates of the data set.
 
 \begin{lemma} There exists an optimal solution $\bf{u^*},\bf{p^*}$ to Problem~\eqref{eq:ERM_2} such that, for all $k= 1, \ldots, K$, $p^*_k \in \{ \hat{p}_i: 1 \leq i\leq I+1 \}$.
 \end{lemma}

\begin{proof}
Let $\hat{p}_i$ be the smallest $p$-coordinate of a data point larger than or equal to $p^*_k$. Replacing $p^*_k$ by $\hat{p}_i $ does not change the value of the objective function.
 \end{proof}

The following proposition shows that we can also   restrict the set of optimal step sizes.
 For a block limited by  $\hat{p}_i$ and $\hat{p}_j$ such that $i < j$, $AV(\hat{p}_i, \hat{p}_j)$,  represents the average of the $x$-coordinates of the data points belonging to that block, i.e.,   $AV(\hat{p}_i, \hat{p}_j) = (\sum_{i \leq h < j} \hat{x}_h)/(j-i)$. Further, let $ER(\hat{p}_i, \hat{p}_j) =\sum_{i \leq h < j} (\hat{x}_h - (AV(\hat{p}_i, \hat{p}_j))^2$.

\begin{prop}
An optimal  solution to Problem~\eqref{eq:ERM_2} is constituted of at most $K'$ blocks, $K' \leq K$, with breakpoints $p^*_k \in  \{ \hat{p}_i: 1 \leq i\leq I+1 \}$ and  step values $u^*_k$ such that:
\begin{enumerate}
\item $u^*_k > u^*_{k+1}$, for $k = 1, \ldots,K'-1$,
\item $u^*_k = AV(p^*_k,p^*_{k+1})$,
\item its objective value for \eqref{eq:OBJ} is equal to $\sum_{k=1}^{K'}ER(p_k^*, p_{k+1}^*)$.
\end{enumerate}

\end{prop} 

\begin{proof}

Consider an optimal solution of problem (\ref{eq:ERM_2})  represented by breakpoints $\{p_k^*\}_{k=1}^{K}$ and step values $\{u_k^*\}_{k=1}^{K}$. Each time that two consecutive blocks, say $k'$ and $k'+1$, have the same function value, i.e., $u_{k'}^{*} =u_{k'+1}^*$,  we can merge them and reduce the number of blocks. Consequently, this optimal solution can be described by $K'$ blocks, with $K' \leq K$, such that $u_k^* > u_{k+1}^*$, for $k = 1, \dots, K' -1$. Given that this solution is globally optimal, it must be locally optimal too, i.e., $u_k^*$ must minimize the contribution of block $k$ to (\ref{eq:OBJ}). In other words,
\begin{align} \label{eq:optimal_step}
u_k^* &\in \arg \min_{u_k} \{\sum_{i: p_k^* \leq \hat{p}_i < p_{k+1}^*} (\hat{x}_i - u_k)^2 : u_{k+1}^* \leq u_k \leq u_{k-1}^*\}  \notag\\ 
u_k^* &\in \{ AV(p_k^*, p_{k+1}^*), u_{k-1}^*, u_{k+1}^*\},
\end{align}

\noindent where $AV(p_k^*, p_{k+1}^*)$ represents the average value of the coordinates $\hat{x}_i$ of data points such that $p_k^* \leq \hat{p}_i < p_{k+1}^*$.
Given that the step values  are all different, it follows that $u_k^* = AV(p_k^*, p_{k+1}^*)$. 

Further,  the contribution of block $k$ to the total error given by \eqref{eq:OBJ} is equal to 
\begin{align}
ER(p_k^*, p_{k+1}^*) =\sum_{i: p_k^* \leq \hat{p}_i < p_{k+1}^*} (\hat{x}_i - (AV(p_k^*, p_{k+1}^*))^2
\end{align}
\end{proof}

We remark that a similar reasoning applies if we consider the least absolute error (that is, the minimization of the sum of absolute values of errors), instead of the least squares. In that case, it suffices to replace the average value of the $\hat{x}_i$-coordinates of the data points such that $p_k^* \leq \hat{p}_i < p_{k+1}^*$ with their \emph{median}.

The above two properties allow to translate problem (\ref{eq:ERM_2}) into a shortest path problem with resource constraints on a particular directed graph $G = (V,A)$  with vertex set $V= \{v_i: 1 \leq i\leq I+1\}$ and edge set $A= \{(v_i, v_j): 1\leq i<j\leq I+1\}$.  

\begin{corollary}
Problem~\eqref{eq:ERM_2} is equivalent to finding a minimum cost path  from $v_1$ to $v_{I+1}$ in graph $G$ with two types of resource constraints:
\begin{enumerate}
    \item the number of arcs in the path is at most $K$,
    \item for any two consecutive arcs $(v_i, v_j)$ and $(v_j, v_h)$ in the path, $AV(\hat{p}_i, \hat{p}_j)> AV(\hat{p}_j, \hat{p}_h)$.
\end{enumerate}
\end{corollary}
\begin{proof}

Each arc $(v_i, v_j)$ corresponds to a block  $[\hat{p}_i, \hat{p}_j)$ with step value $AV(\hat{p}_i, \hat{p}_j)$ and cost $c(v_i,v_j)= ER(\hat{p}_i,\hat{p}_j)$.
Hence, there is a one-to-one correspondence between the paths in $G$ from $v_1$ to $v_{I+1}$ and the set of stepwise functions with breakpoints in  $\{ \hat{p}_i: 1 \leq i\leq I+1 \}$. Adding the two conditions of the corollary ensures that the function is decreasing and contains at most $K$ blocks.

\end{proof}

\subsection{Dynamic programming solution approach }
As previously mentioned, to obtain feasible solutions to LSP, we must impose  two resource constraints. The first one  consists in setting an upper bound $K$ on the numbers of arcs of a path and the second one   excludes  the presence of  consecutive arcs with increasing step values in a path.

The standard approach for solving such a problem consists in using dymamic programming to construct the path from $v_1$ to $v_{I+1}$ progressively, see e.g. \citet{feillet2004exact}. The procedure contains $I$ iterations and at  iteration $i$, partial paths ending in  vertex $v_i \in V'$ are extended by adding one  arc $(v_i,v_j)$.

Further, a  label $l(\pi)$ is associated to each  feasible  partial path $\pi$ from $v_1$  to $v_i \in V'$    specifying the consumption of  the resources.  Here the label  is a triplet $l(\pi) = (c(\pi),k(\pi),st(\pi))$, where $c(\pi)$ denotes the total error of the partial path, $k(\pi)$ its number of arcs and $st(\pi)$ the step value of the block  corresponding to the last  arc of the partial path. On the one hand, the label allows to check whether extending a path $\pi$ ending in $v_i$ by an arc $(v_i, v_j)$ is feasible since we need that $k(\pi) \leq K$ and $st(\pi)> AV(\hat{p}_i, \hat{p}_j)$.
On the other hand, \emph{dominance} between partial paths ending in a same vertex can be exploited.

\begin{definition}
Given two  partial paths $\pi$  and $\pi'$, both ending in  $v_j$, $\pi$ dominates $\pi'$ if $c(\pi)\leq c(\pi'), k(\pi) \leq k(\pi')$ and $st(\pi) \geq st(\pi')$, with at least one strict inequality.
\end{definition}

 If  path $\pi'$ is dominated by some other path, it cannot be part of a  feasible path from $v_1$ to $v_{I+1}$ that  has a strictly better total error. In consequence, all along the execution of the algorithm, we only  need to consider  the partial paths with different and non-dominated  labels. This implies that LSP can be solved in polynomial time.

 \begin{prop}
 An optimal solution to LSP can be found in $\mathcal{O}(KI^3)$ time by solving it as a resource-constrained shortest path problem over graph $G$.
 \end{prop}
 
 \begin{proof}
 Each vertex $v_i$ can be reached by a partial path with  at most $K$ arcs and  the last of these arcs can be associated with at most $i-1$ different step values, namely, $AV(\hat{p}_j, \hat{p}_i)$, with $j = 1, 2, \ldots, i-1$.  Hence, the number of different non-dominated partial paths ending in vertex $v_i$ is in $\mathcal{O}(KI)$. Besides, the number of arcs with $v_i$ as the origin vertex is in $\mathcal{O}(I)$. Consequently, the number of new candidate partial paths generated at iteration $i$ is in $\mathcal{O}(KI^2)$ and the overall complexity of the algorithm is thus $\mathcal{O}(KI^3)$. 
 \end{proof}

Conveniently, we may also take advantage of upper and lower bounds to accelerate the search for the optimal path. Indeed, let $INC$ be the value of a feasible solution to LSP obtained either in some previous iteration of the algorithm or by some other means.  Consider a partial path $\pi$ ending in $v_i$ and let $LB(v_i, k(\pi),st(\pi))$ be a lower bound on the cost of a partial path from $v_i$ to $v_{I+1}$ with at most $K-k(\pi)$ arcs, non-increasing step values and smaller than $st(\pi)$. If $c(\pi) + LB(v_i,k(\pi),st(\pi)) \geq INC$, then the partial path $\pi$ can be directly discarded.
 
A quick valid lower bound can be obtained by relaxing either the condition on the maximum number of  arcs or the constraint on the monotonicity of the step values. In the former case, the problem boils down to an isotonic regression problem on the data points $\{(\hat{p}_j,\hat{x}_j)\}_{j=i}^{I}$. In the latter, it is a lighter shortest path problem from $v_i$ to $v_{I+1}$ in $G$ with at most $K-k(\pi)$ arcs.

The whole procedure is described in Algorithm \ref{algo:shortest_LSP NEW} in which $N_i$ represents the set of labels in the form $ \ell(\pi) = (c(\pi), k(\pi), st(\pi))$ of different and non-dominated partial paths ending in $v_i$ and UB is the initial upper-bound value determined as explained in Section \ref{AS}. Further,  $PRED(\pi)$ is used to store the ``predecessor'' of $\pi$, which is the partial path, say $\pi'$,  that has been extended by one arc to obtain $\pi$. Each time that a new path from $v_{1}$ to $v_{I+1}$ with an objective lower than that of the incumbent solution is found, the variable $OPTIMAL-PATH$ is updated by storing the predecessor of the last node $v_{I+1}$. Once the algorithm terminates, this information allows us to reconstruct the optimal path backwards from $v_{I+1}$ to $v_1$ by starting with $PRED(OPTIMAL-PATH)$. This way, the optimal solution to LSP is eventually retrieved.

\begin{algorithm}[htpb]
\caption{Shortest path algorithm for LSP}
\label{algo:shortest_LSP NEW}
\begin{algorithmic}[1]
\STATE { \bf Initialization:} $N_1 = \{(0,0,\max_{i \in I} \hat{x}_i+1)\}$, $INC = UB$
\FOR{$i \in \{1,\ldots, I \}$}

	\WHILE{$N_i \neq \emptyset$}
		\STATE Select $\pi^* \in \mbox{arg}\min_{\ell(\pi)\in N_i} \{ c(\pi)  \}$ and remove $\ell(\pi^*)$ from $N_i$
		\IF{$c(\pi^*)> INC$}
			\STATE $N_i = \emptyset$
		\ELSIF{ $c(\pi^*) + c(v_i, v_{I+1}) < INC$ \AND $AV(p_i,p_{I+1}) < st(\pi^*)$}
				\STATE $PRED(OPTIMAL-PATH) = \pi^*$, $INC = c(\pi^*) + c(v_i, v_{I+1})$
		\ENDIF
		\IF{$k(\pi^*) < K-1$}
	        \FOR{$h \in i+1, \ldots, I$}
				\IF{$AV(p_i,p_h) < st(\pi^*)$ \AND $c(\pi^*) +c(v_i, v_h) +  LB(v_h,k(\pi^*)+1, AV(p_i,p_h)) < INC$}
					\STATE $new = (c(\pi^*) +c(v_i, v_h), k(\pi^*)+1, AV(p_i,p_h))$, $PRED(new) = \pi^*$ 
					\IF{$new \notin N_h$}
						\STATE Add label $new$ to $N_h$ if it is not dominated.
						\STATE Delete all dominated labels.
					\ENDIF
				\ENDIF
			\ENDFOR
			\ENDIF
		
	\ENDWHILE	
\ENDFOR
\end{algorithmic}
\end{algorithm}

\section{Acceleration strategies}\label{AS}

Despite the fact that LSP can be solved in polynomial time, computing the optimal solution can be expensive for realistic instances. The overall solution time relies heavily on how tight  the upper bound $INC$ and the lower bounds $LB(\cdot)$ are. In this section we discuss strategies to find good bounds that are easy to compute.

\subsection{Computing an upper bound: Combining isotonic regression with adjacency-constrained data clustering}\label{sec:UB}

Feasible solutions for problem~\eqref{eq:ERM} provide us with an upper bound on the optimal error that can help us reduce the computational burden of the shortest path problem presented in Section \ref{s:shortest_algorithm}. One efficient procedure to compute a tight upper bound runs as follows:
\begin{enumerate}
    \item We use isotonic regression to fit a monotone stepwise function to the original data set. However, one should expect the number of blocks of this fit to be higher than $K$.
    \item We reduce the number of blocks of the output of the isotonic regression to $K$ by grouping the consumption values of the isotonic fit into $K$ clusters. For this purpose, we use the fast greedy algorithm proposed in~\citet{pineda2018chronological} for adjacency-constrained hierarchical clustering.
    \item The step value is computed as the average consumption of the isotonic fit values within each of the $K$ clusters obtained in the previous point.
\end{enumerate}

\begin{figure}[htpb]
    \centering
    \includegraphics[width = \textwidth]{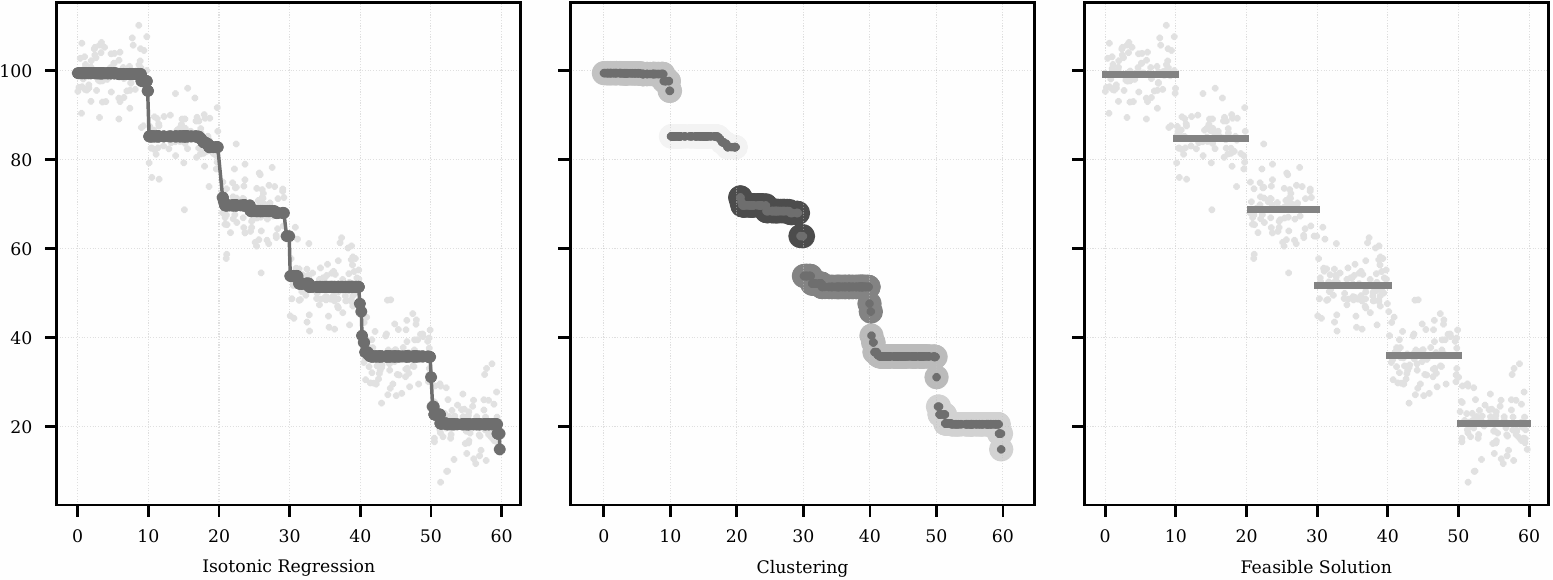}
    \caption{Algorithm to compute a feasible solution and an upper bound to the minimum error.}
    \label{fig:upper}
\end{figure}

The procedure above yields a monotone stepwise function with $K$ pieces that is a feasible solution to LSP. This methodology is depicted in Figure \ref{fig:upper}, with each subfigure representing one of the actions described above, from left to the right. We implement the calculation of the so-obtained upper bound on Python, using the isotonic regression and the agglomerative clustering functions of package \emph{Scikit-learn}, see \citet{scikit-learn}.

\subsection{Computing a lower bound}\label{CLB}

Lower bounds are useful in several ways. First, as we discussed in Section~\ref{s:shortest_algorithm}, they prevent the shortest-path algorithm from creating sub-optimal labels. To do so, it is necessary to compute lower bounds for each partial path $\pi$. Depending on the method, this can be computationally expensive. Second, lower bounds give a guarantee of how far any feasible solution is from the optimal one. We obtain these lower bounds by relaxing either the  constraint on the number of blocks/arcs or the monotonicity constraint of the fitted curve. 

\paragraph{Relaxing the  constraint on the number of arcs in the path: The isotonic fitted curve}

When  the number of blocks is not limited, problem  LSP  is equivalent to the well-known isotonic fit, (\citet{fielding1974statistical}). Isotonic regression can be solved in linear time (\citet{best1990active}). Given the efficiency of this method, we generate lower bounds for any partial path $\pi$ by computing one lower bound for each vertex $v_i$. In other words, we calculate $LB(v_i, k(\pi), st(\pi))$ as $LB(v_i)$ for each partial path $\pi$. The total time to compute this lower  bound for all  $v_i$ is  in $O(I^2)$. We use the isotonic regression function implemented in the Python package \emph{Scikit-learn} to this end (\citet{scikit-learn}).

\paragraph{Relaxing the monotonicity constraint} As mentioned in Section \ref{s:shortest_algorithm}, another lower bound can be obtained by relaxing the monotonicity constraint on the step values. Then,  for a given partial path $\pi$,  a lower bound can be computed by solving a shortest path problem with at most $K-k(\pi)$ arcs from $v_i$ to $v_{I+1}$  in $G'$.  One can determine the fitting error associated with  the shortest paths  from all vertices $v_i$ to the sink $v_{I+1}$  and containing at most $k$ arcs, for all $k= 1, \ldots,K$, with dynamic programming. The corresponding  total computing time is  in $ O(K\vert A \vert)=O(KI^2)$ (essentially, if we disregard the monotonicity constraint, our problem translates into a standard cardinality-constrained shortest problem whose computational complexity is known to be in $\mathcal{O}(K|A|)$, with $|A| = \mathcal{O}(I^2)$ in our case). 

In terms of implementation, to relax the monotonicity constraint is equivalent to suppressing the third component \emph{st} of each label and the corresponding monotonicity conditions in line 7 and 12 in Algorithm \ref{algo:shortest_LSP NEW}.
In particular, for the path that corresponds to the initial label $(0,0)$ and contains the single vertex $v_1$,  this shortest path problem  returns a lower bound on the minimum fitting error. In that case, besides, if the resulting function turns out to be non-increasing, then it must be optimal to LSP. 

We end this subsection with a remark on the so modified algorithm: The minimum cost computed over all the partial paths reaching a layer $i$ in the modified algorithm is a lower bound on that very same cost in the original Algorithm 1. Consequently, we can get an even tighter lower bound by running Algorithm 1 with the monotonicity constraint dropped and with $LB(v_h,k(\pi^*)+1, AV(p_i,p_h))$ in line 12 of the pseudocode given by the isotonic fit. The total cost at termination does not necessarily correspond to the fitting error of the optimal, possibly non-monotone, stepwise curve, but it is a still valid lower bound on the cumulative error of LSP. 
This is indeed the lower bound (obtained from dropping the monotonicity constraint on the optimal fitted curve) that we will consider in the numerical experiments below.

\subsection{Imposing constraints on the length of steps}

In some cases, it may be interesting and practically useful to impose constraints on the length of the function blocks, i.e., to restrict the set of functions $\mathcal F$ to decreasing stepwise functions with a step length bigger than \texttt{step\_min}. This is equivalent to adding the following set of constraints to model \eqref{eq:ERM_2}:
\begin{equation}
   \mbox{\texttt{step\_min}}   \le p_{k+1} - p_{k}, \quad \forall k \le K \label{eq:step_constraint}
\end{equation}

Moreover, Algorithm \ref{algo:shortest_LSP NEW} can still be used after removing from the arc set $A$ all   $(v_i,v_j)$ for which the corresponding $p$-coordinates violate condition (\ref{eq:step_constraint}). As a result, the number of operations to compute the optimal solution to LSP decreases. We show  the impact of imposing this type of constraint experimentally in the next section.

We remark that the upper bound described in Section~\ref{sec:UB} may no longer be feasible after enforcing the constraint on the minimum step length. Nevertheless, if that is the case, we can always gradually decrease $K$ in the algorithm outlined in that section until a valid upper bound is eventually recovered.

\section{Numerical experiments}\label{NE}

Next we run a series of numerical experiments to test the effectiveness and performance of the proposed algorithm under different settings. To this end, we first use synthetic data sets to assess the sensitivity of the algorithm performance to the maximum number of steps, the noise level in the input data and the sample size. Subsequently, we consider a realistic data set consisting of price-power measurements at the main substation of a distribution power grid that includes distributed energy resources. This data can be download from \citep{OASYS2020segisoreg}. 

All the numerical experiments have been conducted on a laptop equipped  with  a  Intel  Core i5-4200 CPU processor, with 2.80 gigahertz 2-core, and 8 gigabytes of RAM memory.  The operating system is 64-bit Windows 8.1.  Codes were implemented in Python 3.8.

\subsection{Synthetic data sets}

\begin{figure}
    \centering
\begin{subfigure}[t]{0.48\textwidth}
\includegraphics[width=\textwidth]{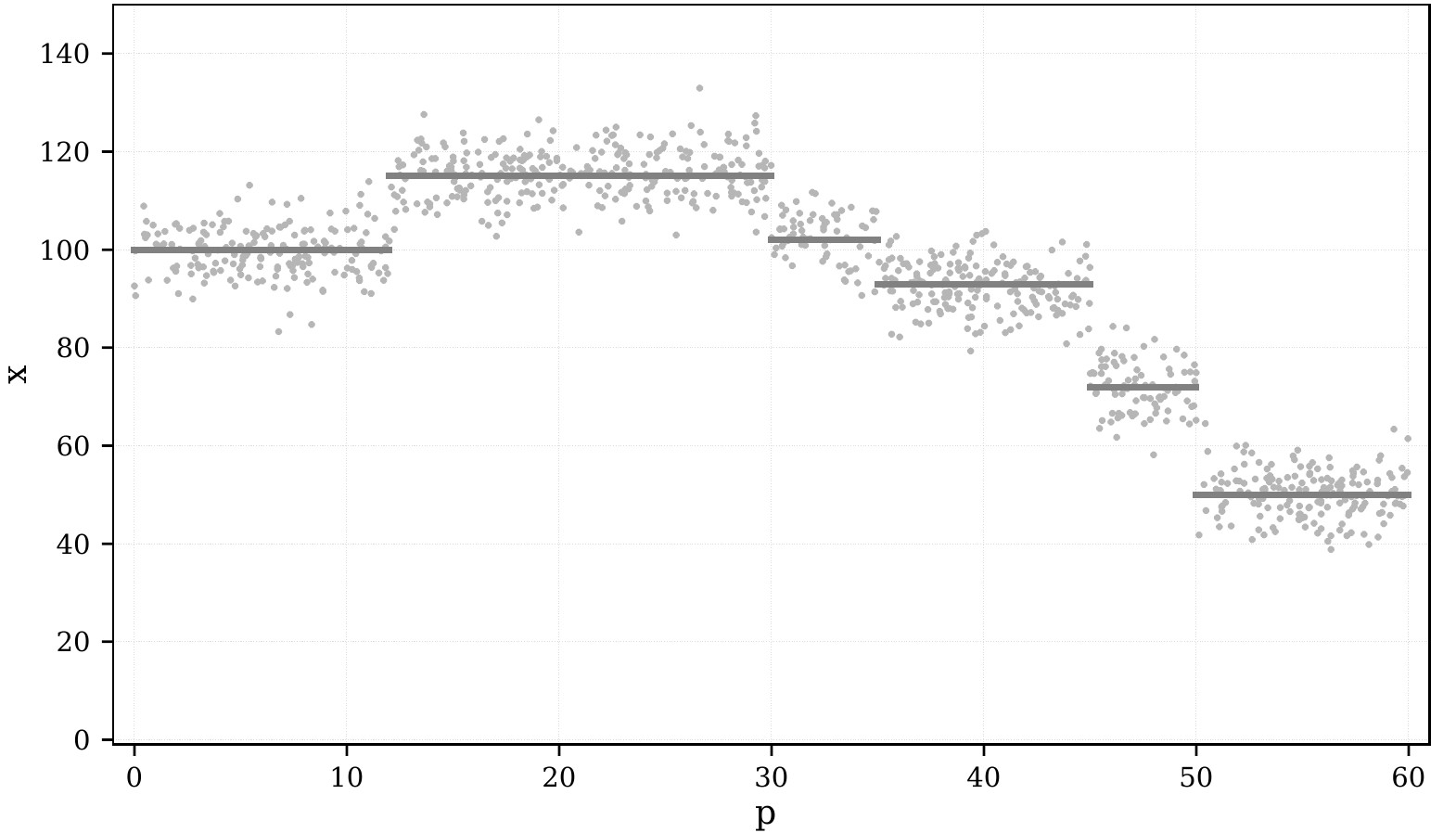}
\caption{$\sigma = 5$}
\label{fig:noise_5}
\end{subfigure}
\begin{subfigure}[t]{0.48\textwidth}
\includegraphics[width=\textwidth]{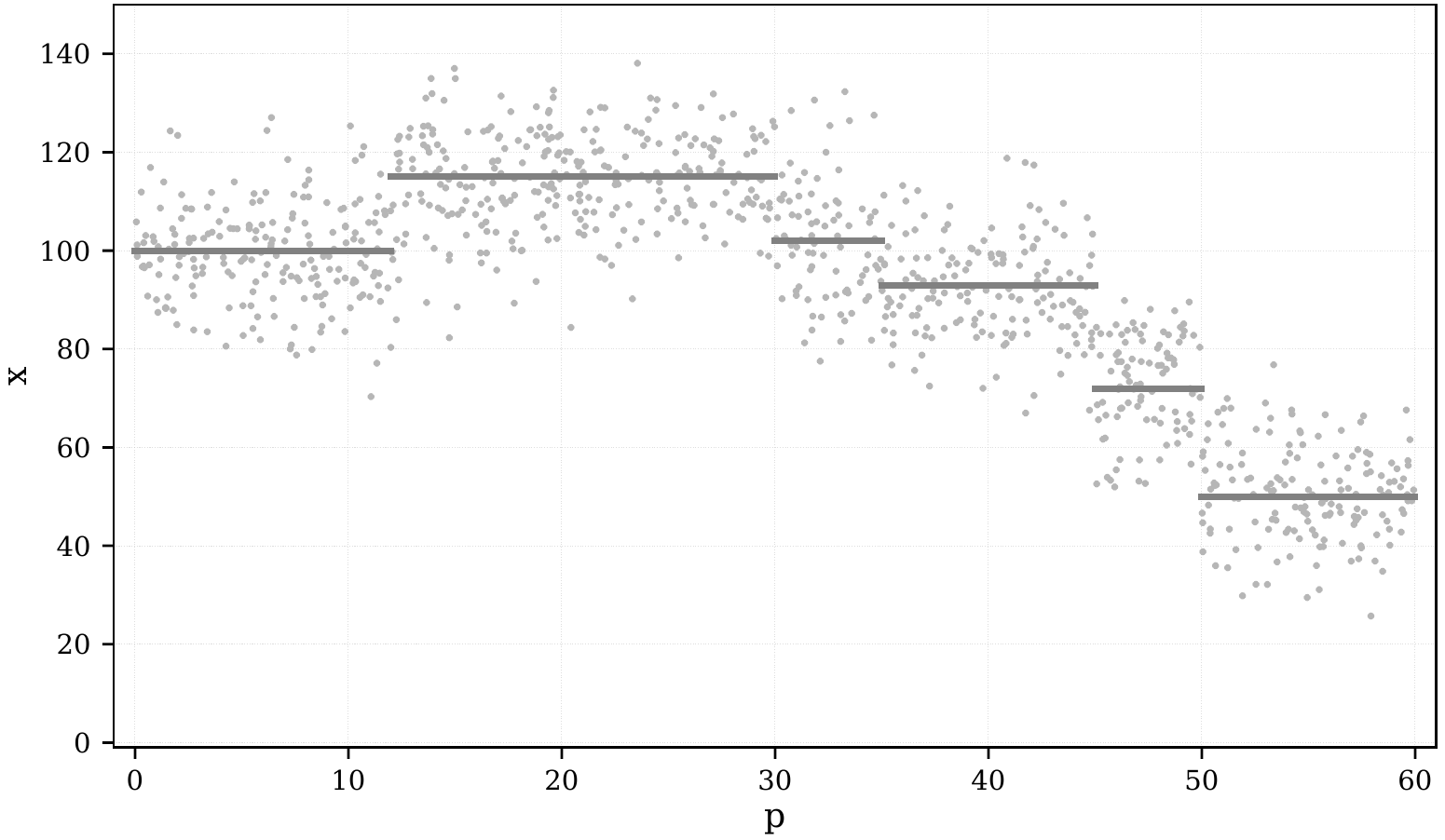}
\caption{$\sigma = 10$}
\label{fig:noise_10}
\end{subfigure}
 \caption{Synthetic data for different noise levels.}
        \label{fig:noises}
\end{figure}

We first test our algorithm and the effectiveness of the acceleration strategies described above on a controlled experiment, where we know the true data-generating distribution. More specifically, the response variable $x$ is given by
\begin{equation}
x = f^*(p)+ \varepsilon \label{numerics: true model}
\end{equation}
where $\varepsilon$ is a Gaussian noise of zero mean and standard deviation $\sigma$ and $f^*$ is the stepwise function depicted in Figure~\ref{fig:noises}, that is,
\begin{equation}
f^*(p) = \sum_{k=1}^{6}{u_k^* \mathbb{I}_{[p^*_k,p^*_{k+1})}(p)} \label{eq:true_response}
\end{equation}
with $u_k^*$ and $[p^*_k,p^*_{k+1})$, $k=1,\ldots, 6$, provided in Table~\ref{tab:synthetic_curve}.
\begin{table}
    \centering
    \begin{tabular}{ccc}
    \hline
    $k$-$th$ block  & $[p^*_k,p^*_{k+1})$  & $u_k^*$  \\
    \hline
1  & [0,12) & 100 \\
2  & [12,30) & 115\\
3  & [30,35) & 102\\
4 & [35,45) & 93\\
5 & [45,50) & 72\\
6 & [50,60] & 50\\
    \hline
    \end{tabular}
    \caption{Stepwise characterization of the \emph{true} relationship between the response $x$ and the covariate $p$.} 
    \label{tab:synthetic_curve}
\end{table}
Notice that $f^*$ is a stepwise function made up of six blocks of different sizes. Furthermore, $f^*$ is neither increasing, nor decreasing in its entire domain. For illustration purposes, Figures \ref{fig:noise_5} and \ref{fig:noise_10} plot 1000 data points $\{(\hat{p}_i,\hat{x}_i)\}_{i=1}^{I}$ randomly generated for two different noise levels, namely, $\sigma=5$ and $\sigma=10$, respectively. Both figures also include the true function \eqref{eq:true_response} to compute the response variable $x$. 

Next, we run our shortest path algorithm using datasets of 1000 points that are randomly generated from \eqref{numerics: true model} for a noise level $\sigma$ taking the values of the natural numbers between zero and ten. Besides, the number of arcs $K$ is set to six, which is the true number of blocks of the function that relates the response variable $x$ to the covariate $p$. Results for all these cases are collated in Table \ref{tab:impact_noise} and include the aggregated square error (Error) and three different computational times (with a maximum value of 12 hours):
\begin{itemize}
    \item[-] T$^{\rm ISO}$: Computational time of  the proposed shortest path algorithm including the upper bound  discussed in Section \ref{sec:UB} and the lower bound per layer provided by the isotonic fit.
    \item[-] T$^{\rm RLX}$: This computational time is obtained as follows. Let T$^{\rm W/O}$ denote the time needed to run the proposed shortest path algorithm \emph{without the monotonicity constraint}, but including the upper bound  discussed in Section \ref{sec:UB} and the lower bounds (one per layer) provided by the isotonic regression fits. As mentioned in Section~\ref{CLB}, the cost of this shortest path constitutes a valid lower bound on the cumulative fitting error associated with the optimal monotone curve. Actually, if this shortest path leads to a nonincreasing curve, then this is the optimal one. In this case, we set T$^{\rm RLX}$ = T$^{\rm W/O}$ (because the optimal fit has been found). Otherwise, we need to rerun our algorithm \emph{with the monotonicity constraint} back in force, and therefore, we set T$^{\rm RLX}$ = T$^{\rm W/O}$ + T$^{\rm ISO}$. 
    \item[-] T$^{\rm NOB}$: Computational time of the proposed shortest path algorithm if no acceleration strategies are employed. 
\end{itemize}

\begin{table}[htpb]
    \centering
    \begin{tabular}{ccccccc}
    \hline
    $I$ & $\sigma$ & $K$ & Error & T$^{\rm ISO}$(s) & T$^{\rm RLX}$(s) & T$^{\rm NOB}$(s) \\
    \hline
1000 & 0 & 6 &  28081	& 110 & 164 & 13053 \\
1000 & 1 & 6 &  27491	& 647 & 837 & $>$43200     \\
1000 & 2 & 6 &  30716	& 1296 & 1569 & 39205 \\
1000 & 3 & 6 &  35320	& 441 & 644 & 33675 \\
1000 & 4 & 6 &  41852	& 1258 & 1501 & 41128\\
1000 & 5 & 6 &  52919	& 1038 & 1177 & 37186 \\
1000 & 6 & 6 &  64416	& 1400 & 1547 & $>$43200 \\
1000 & 7 & 6 &  67985	& 1222 & 1447 & $>$43200\\
1000 & 8 & 6 &  102860	& 1998 & 2346 & $>$43200 \\
1000 & 9 & 6 &  113847	& 1570 & 1903 & 42156 \\
1000 & 10 & 6 & 123856	& 1162 & 1143 & 32689 \\
    \hline
    \end{tabular}
    \caption{Impact of data noise}
    \label{tab:impact_noise}
\end{table}

By comparing the computational times T$^{\rm NOB}$ and T$^{\rm ISO}$ of Table~\ref{tab:impact_noise}, we can conclude that the use of the proposed upper and lower bounds has a tremendous impact on the ability of the algorithm to quickly identify the globally optimal curve to be fitted. Besides, these results also reveal that our algorithm is robust to the level of noise, since the computational time T$^{\rm ISO}$ is relatively stable as noise increases. Finally, the lower bound provided by the relaxation of the monotonicity constraint is, nevertheless, of little value for this instance, in which T$^{\rm RLX}$ is higher than T$^{\rm ISO}$ for most noise levels. We will see, however, that this lower bound can be useful when the data features a sufficiently marked monotonic layout.

In order to better understand the intuition behind the acceleration strategies described in Section \ref{s:shortest_algorithm} and their impact on the computational time of the shortest path problem proposed in Section \ref{AS}, Figure ~\ref{fig:bounds} displays the following:
\begin{figure}[]
    \centering
    \includegraphics[width = \textwidth]{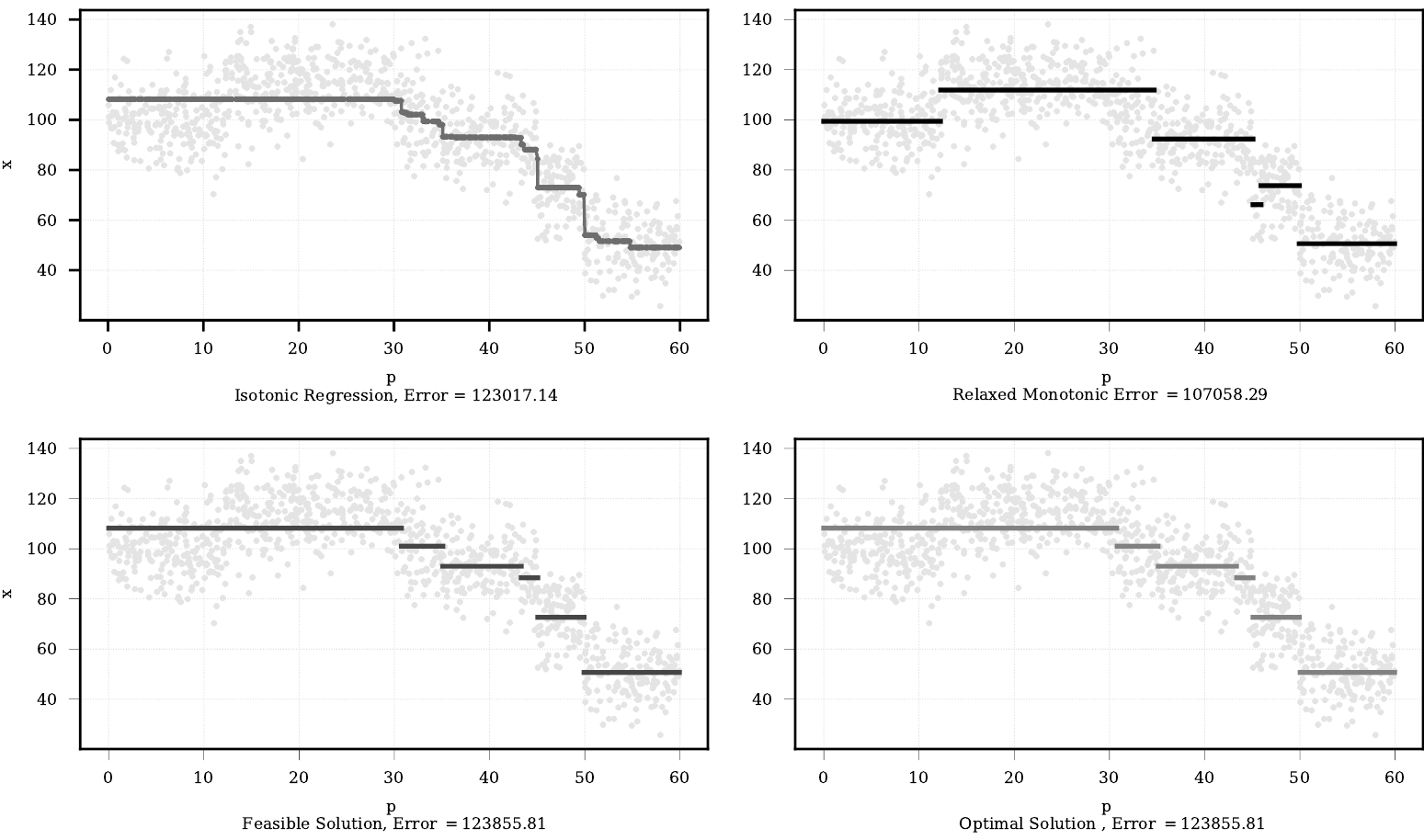}
    \caption{Illustration of the various bounds and the optimal solution for the synthetic data with noise $\sigma =10$ and $K=6$}
    \label{fig:bounds}
\end{figure}
\begin{itemize}
    \item[-] Top left plot: The curve provided by the isotonic regression. As observed, the isotonic fit is non-increasing, but the number of steps is higher than $K =6$. Therefore, the aggregated squared error associated with this curve can be used to lower-bound the optimal solution.
    \item[-] Bottom left plot: The curve obtained through the adjacency-constrained hierarchical clustering technique using  the isotonic fit as input. The number of blocks is equal to six and the monotonic condition is also satisfied. Therefore, this curve represents a feasible solution in LSP and its accrued squared error is a valid upper bound of the optimal objective function value.
    \item[-] Top right plot: The curve computed by the relaxed shortest path algorithm without the monotonicity constraint. The number of blocks is also equal to six, but the curve is not monotone. Therefore, the corresponding aggregate squared error can also be used as a lower bound.
    \item[-] Bottom right plot: The global optimal solution obtained by the proposed shortest path algorithm. 
\end{itemize}

Interestingly, the optimal solution of this particular instance coincides with that resulting from the combination of the isotonic regression and the structured hierarchical clustering. Furthermore, the lower bound provided by the isotonic fit is notoriously tight, which is, most likely, the reason behind the good performance exhibited by our algorithm. Notice that the lower bound achieved by relaxing the monotonicity constraint is significantly less tight than the one given by the isotonic regression fit.

\begin{figure}[htpb]
    \centering
    \includegraphics[width = \textwidth]{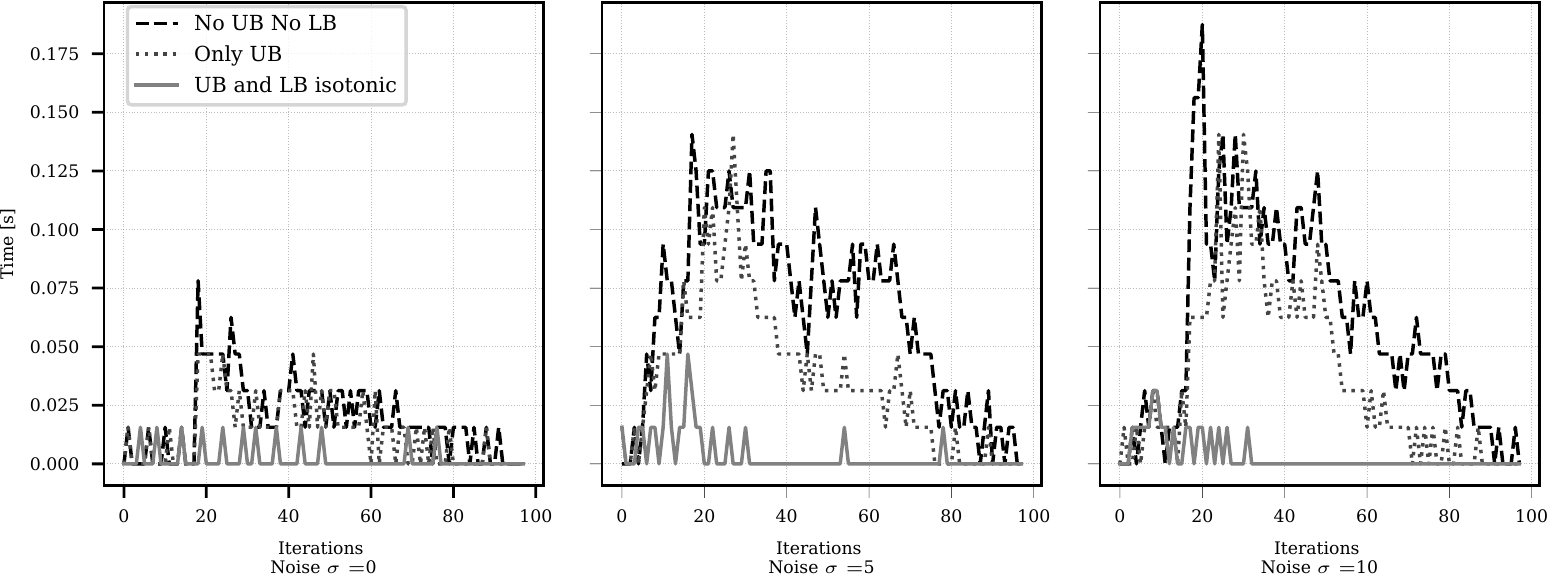}
    \caption{Effect of bounds on computational burden per iteration. Example with 100 data points and $K=6$}
    \label{fig:time_per_iter}
\end{figure}

To further illustrate how the proposed bounds can accelerate the solution of the shortest path algorithm, Figure~\ref{fig:time_per_iter} shows the time spent per iteration by our algorithm for noise levels $\sigma = 0$, $5$ and $10$. In this figure, the dashed plots refer to the raw implementation of the algorithm, i.e., with no bounds; the  dotted lines correspond to the version of the algorithm where only the proposed upper bound is used; finally, the solid plots provide the time our algorithm spends per iteration when both the lower and the upper bounds are exploited.
%
It is apparent that using bounds in the proposed methodology has a remarkable beneficial effect on the algorithm performance, to such an extent that the joint use of both bounds  manages to immunize the algorithm against the noise. Indeed, our upper and lower bounds noticeably reduce the number of labels that Algorithm 1 generates in the intermediate layers of the graph. In the limiting case where there is no noise, both the upper and the lower bounds coincide with the optimal fit and no intermediate label is generated at all, thus taking a marginal amount of time per iteration. As the noise level is increased, more and more intermediate labels are to be handled, which essentially tells us that the optimization underlying the regression problem becomes harder and harder to perform. Furthermore, as can be inferred from the plots of Figure~\ref{fig:time_per_iter}, 
the inclusion of the upper bound only is not enough to keep the computational burden of our algorithm per iteration low, because of the high amount of labels that are produced in the first layers of the graph. It is the synergistic effect of the lower and upper bounds which prevents the number of labels in the early stages of our algorithm from exploding.  

In what follows, we omit computational time $T^{\rm NOB}$, since it has become clear that our algorithm runs much faster when combined with the proposed acceleration strategies.

Now we fix the noise level $\sigma$ to five and change the sample size instead. Still we have $K = 6$. The comparison results of T$^{\rm ISO}$ and T$^{\rm RLX}$ for this new experiment are collated in Table \ref{tab:impact_size}. Naturally, the aggregate squared error and the solution times increase with the sample size $I$. However, our algorithm appears to scale relatively well, given its theoretical complexity. 

\begin{table}[htpb]
    \centering
    \begin{tabular}{cccccc}
    \hline
    $I$ & $\sigma$ & $K$ &Error & T$^{\rm ISO}$(s) & T$^{\rm RLX}$(s) \\
    \hline
100  & 5 & 6 &  5074	&0&0 \\
200  & 5 & 6 &  11176	&2&3\\
500  & 5 & 6 &  28810	&73&99\\
1000 & 5 & 6 &  52919	& 1038&1177\\
2000 & 5 & 6 &  112410	&9092&10726\\
    \hline
    \end{tabular}
    \caption{Impact of sample size}
    \label{tab:impact_size}
\end{table}

Finally, we fix the sample size to 1000 and the noise level of the data to five, and change the maximum number of arcs $K$ our algorithm may use to reduce the error. The so obtained results are compiled in Table \ref{tab:impact_blocks}. As expected, by increasing the maximum number of arcs $K$ (also referred to as \emph{number of blocks}), we enrich the family $\mathcal{F}$ of non-increasing stepwise functions we consider and thus, the error of the data fitting is reduced. If the number of blocks $K$ is lower than or equal to four, the solution obtained by the relaxed shortest path without the monotonicity constraint happens to be non-incresing and thus, optimal. This explains why T$^{\rm RLX}$ is significantly lower than T$^{\rm ISO}$ if $ K\leq 4$. On the contrary, if $K$ is higher or equal to five, relaxing the monotonicity constraint leads to non-monotone solutions in order to adapt as much as possible to the original function, which is also non-monotone. In such cases, the time required to compute the lower bound through the relaxed shortest path problem is significantly higher than the time savings originated by such lower bound and consequently, T$^{\rm ISO}$ is lower than T$^{\rm RLX}$ for $ K\geq 5$.

\begin{table}[]
    \centering
    \begin{tabular}{cccccc}
    \hline
    $I$ & $\sigma$ & $K$ & Error & T$^{\rm ISO}$(s) & T$^{\rm RLX}$(s) \\
    \hline
1000 & 5 & 2 &118510	&0&0\\
1000 & 5 & 3 &81549	    &42&26\\
1000 & 5 & 4 &55242	    &585&52\\
1000 & 5 & 5 &53275	    &744&884\\
1000 & 5 & 10 &52528    &1799&1821\\
1000 & 5 & 12 &52493    &1622&1798\\
1000 & 5 & 14 &52479    &1817&2058\\
1000 & 5 & 16 &52471    &1951&2361\\
1000 & 5 & 18 &52466    &2234&2791\\
1000 & 5 & 20 &52465    &2734&3637\\
\hline
    \end{tabular}
    \caption{Impact of the maximum number of arcs $K$}
    \label{tab:impact_blocks}
\end{table}

\subsection{Realistic application: Estimating the bidding curve of a pool of flexible consumers}

\begin{figure}[htpb]
    \centering
    \includegraphics[width = \textwidth]{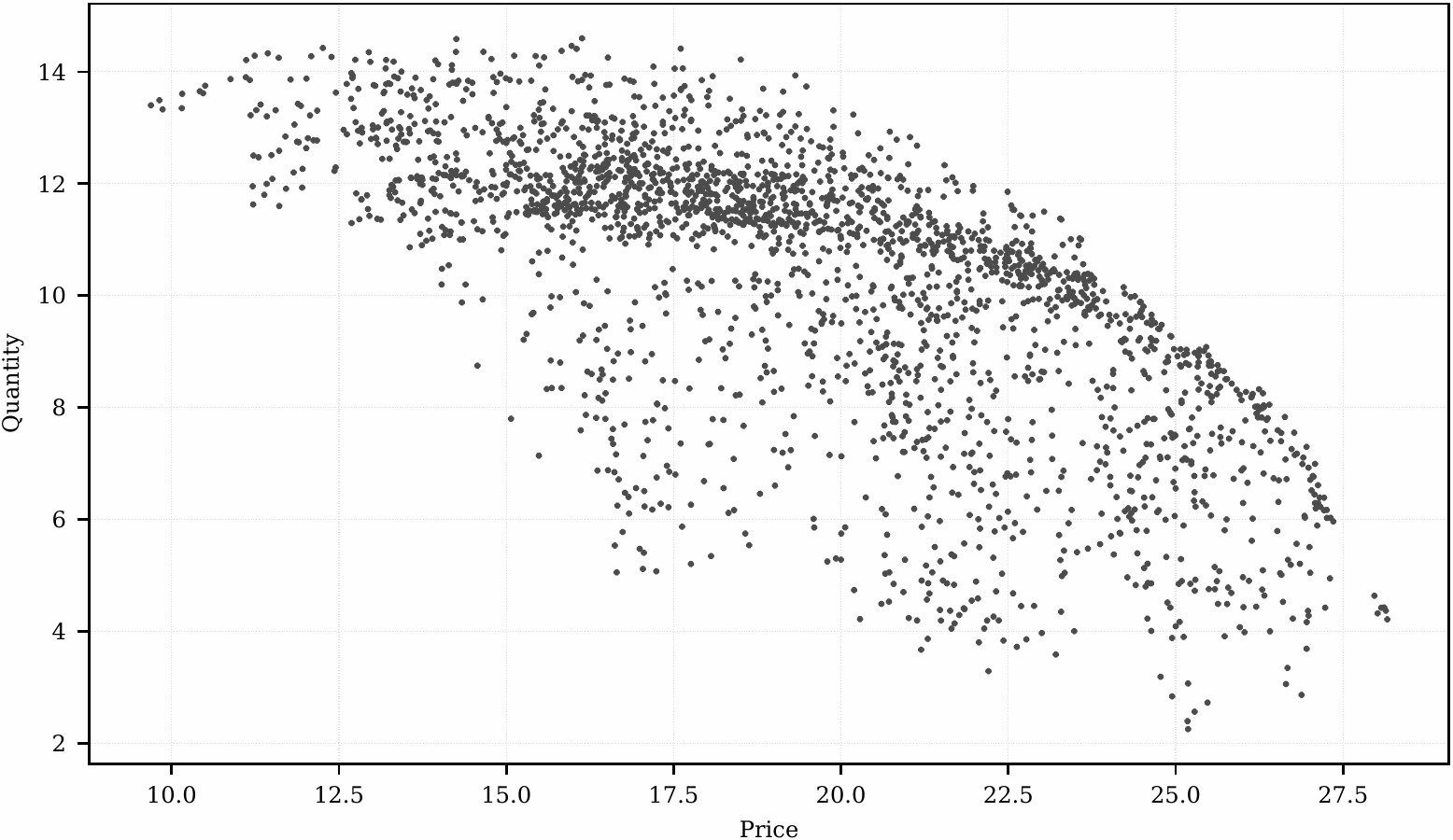}
    \caption{Price-consumption data from realistic application}
    \label{fig:realistic_data}
\end{figure}

Here we consider the problem of estimating the price-response of a cluster of flexible consumers of electricity, that is, how much energy the cluster consumes as a function of the electricity price. Similar instances of this problem has been considered, for example, in~\citet{aswani2018inverse, saez2017short, saez2016data}. In our particular case, these consumers are located within the 33-bus radial distribution grid described in \citet{Hassan2018}. The distribution network includes 8 solar generating units whose power output varies through time according to weather conditions, and 32 flexible consumers able to adapt their consumption to electricity prices as modeled in \citet{Mieth2020}. As proposed in \cite{Mieth2018}, a LinDistFlow modeling approach is used to account for both voltage and line capacity limits. Detailed data about all parameters of the distribution network is available at \citep{OASYS2020segisoreg}. The distribution grid receives a nodal price at the main substation, to which the consumers react according to their energy needs, generation assets, and sensitivity to the electricity cost. The aggregate amount of energy demanded by the pool, paired with the nodal price (at the main substation) that induced such a demand, constitutes an observation and form a data point on the plane. The collection of the 2400 observations at our disposal are plotted in Figure~\ref{fig:realistic_data} and can be downloaded from \citep{OASYS2020segisoreg}. Besides the electricity price, the operation of the distribution network is also affected by other factors such as the varying solar power generation and therefore, similar price signals may yield quite different consumption levels.

According to current rules in many day-ahead electricity markets worldwide, consumers must submit a stepwise and nonincreasing bidding curve indicating their demand levels as a function of the electricity price. Besides, deviations with respect to the declared consumption quantities are to be penalized in the real-time/balancing market. In most electricity markets, over- and under-consumption are equally penalized according to a single-price settlement. Under these conditions, the cluster of consumers is interested in finding the stepwise nonincreasing function that minimizes the mean squared error with respect to the data provided in Figure~\ref{fig:realistic_data}. If positive and negative energy deviations are priced differently according to a two-price balancing settlement, the proposed procedure can also be used by replacing $AV(p_k^*, p_{k+1}^*)$ in \eqref{eq:optimal_step} with the appropriate empirical quantile. For all these reasons, this problem represents a natural practical application of the mathematical problem described in Section \ref{PF}.

The cumulative squared error of the curve fit provided by our algorithm for a different number of arcs (or steps) is compiled in Table~\ref{tab:impact_realistic}. This table also shows the solution times T$^{\rm ISO}$ and T$^{\rm RLX}$ defined in the previous example, the initial upper bound and the lower bound obtained by relaxing the monotonicity constraint from which our algorithm starts to iterate. From these bounds, we can compute the optimality gap ${\rm GAP}_0 = \frac{UB-LB}{LB}100\%$ at the beginning of the algorithm, which we include in such a table too. We omit time T$^{\rm NOB}$, as the raw algorithm is unable to deliver the optimal solution within a day in most cases, which proves the computational efficiency of the proposed acceleration strategies for our shortest path algorithm. As a matter of fact, the initial optimality gap that our algorithm needs to close is always below 0.25\%, which reveals that the heuristic procedures we have devised to construct a (feasible) upper bound and a tight lower bound are remarkably good. Consequently, if the bidding curve is to be determined very frequently to participate in intra-day trading floors, then using the solution provided by the proposed heuristic procedures may be a good compromise between accuracy and computational time. Conversely, using the proposed exact procedure is justified to obtain the bidding curve to be submitted to a day-ahead electricity market once every 24 hours.

For those cases in which Algorithm 1 reaches the maximum time limit and thus, is terminated without having certified that the optimal curve fit has been found, we include, within parentheses, the optimality gap at termination. This optimality gap is calculated from the best upper and lower bounds on the optimal solution that are available after the time limit has expired. In the case of the experiment associated with time T$^{\rm ISO}$, which only considers the lower bounds given by the isotonic fit, the best lower bound at termination is computed as follows:
\begin{enumerate}
    \item Let $i'$ be the layer being processed by the loop \emph{for} in line 2 in the pseudocode of Algorithm 1 at termination. Consider each feasible path reaching layer $i'-1$ and the cost accrued by this path until that layer. Increase this cost by the error of the isotonic fit from layer $i'-1$ to the last one $I+1$. Denote the result as the \emph{extended cost} of a feasible path at layer $i'-1$.
    \item Compute the minimum extended cost for each layer $i \in \{1, 2, \ldots, i'-1\}$.
    \item The best lower bound is then given by the maximum over layers $\{1, 2, \ldots,$ $i'-1\}$ of their associated minimum extended cost. 
    
\end{enumerate}

\begin{table}[htpb]
\small
    \centering
    \begin{tabular}{crrrrrr}
    \hline
    $K$ &Error & T$^{\rm ISO}$(s) & T$^{\rm RLX}$(s) & LB & UB & ${\rm GAP}_0$(\%) \\
    \hline
    1 & 14901 & 1.1 & 1.1 & 14901 & 14901 &  0\\
    2 & 9201 & 2.5 & 2.2 & 9201 & 9201     & 0\\
    3 & 8213 & 731.4 & 581.3 & 8213 & 8222 & 0.11\\
    4 & 7726 & 30602 & 1077.2 & 7726 & 7732 & 0.08 \\
    5 & 7596 & 65995 & 1881.3 & 7596 &  7602 & 0.08\\
    6 & 7502 & $>$86400 (2.71\%) & 2477.3 & 7502 & 7519 & 0.23\\
    7 & 7442 & $>$86400 (1.88\%)  & 2895.8 & 7442 & 7448 & 0.08\\
    8 & - & $>$86400 (1.31\%)& $>$86400 (0.14\%)& 7392  & 7402 & 0.14\\
    \hline
    \end{tabular}
    \caption{Realistic application: Cumulative squared error, solution times, bounds and optimality gaps}
    \label{tab:impact_realistic}
\end{table}

 Once again, the optimality gaps provided within parentheses confirm that the feasible solution we construct at the beginning of the algorithm, by modifying the isotonic fit through adjacency-constrained data clustering, is nearly optimal and that the lower bound given by relaxing the monotonicity constraint is hard to beat for this data set.

Results in Table \ref{tab:impact_realistic} also show that the accrued fitting error decreases with the number of blocks, since the family $\mathcal{F}$ of non-increasing stepwise functions becomes larger as $K$ is augmented. Nonetheless, the reduction in the fitting error we get by increasing $K$ rapidly plateaus after $K > 4$.  In practice,  the number of pieces $K$ should be treated as a hyperparemeter determining the complexity of our statistical model. As such, we can use standard strategies available in the machine-learning literature for hyperparameter tuning to properly set $K$. For instance, we can borrow the popular \emph{elbow criterion} from the realm of data clustering for this purpose. According to this criterion, the explained variation as a function of $K$ is plotted, and the elbow of the curve is taken as a good value for $K$. In the present case, the elbow is clearly placed on the value $K = 4$. Alternatively, we can also use more sophisticated cross-validation procedures to this very same aim \cite[Ch. 7]{hastie2009elements}. In principle, low values of $K$ should be preferred to favor model simplicity and avoid overfitting.

\begin{figure}[htpb]
    \centering
       \includegraphics[width = \textwidth]{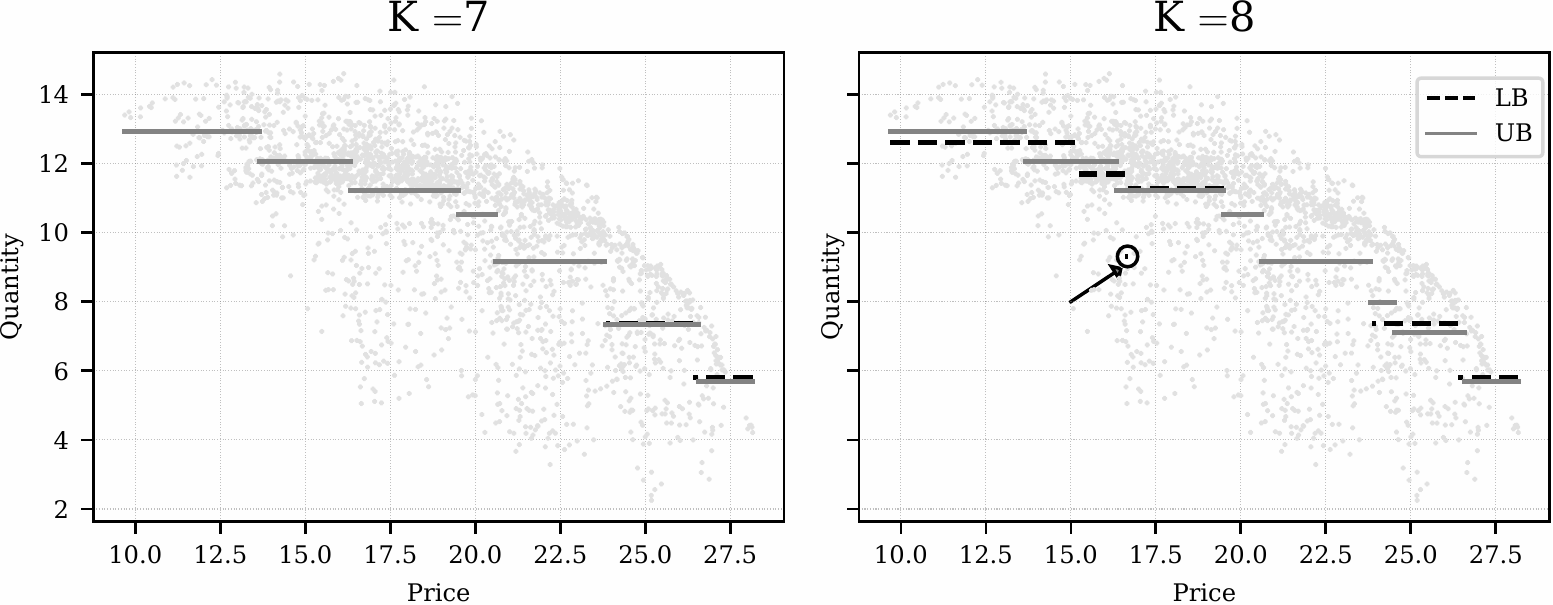}
    \caption{Lower and upper bound solutions for $K=7$ and $K=8$. Notice the tiny step that appears in the fit provided by the lower bound for $K=8$}
    \label{fig:realistic_k78}
\end{figure}

\begin{figure}[htpb]
    \centering
      \includegraphics[width = \textwidth]{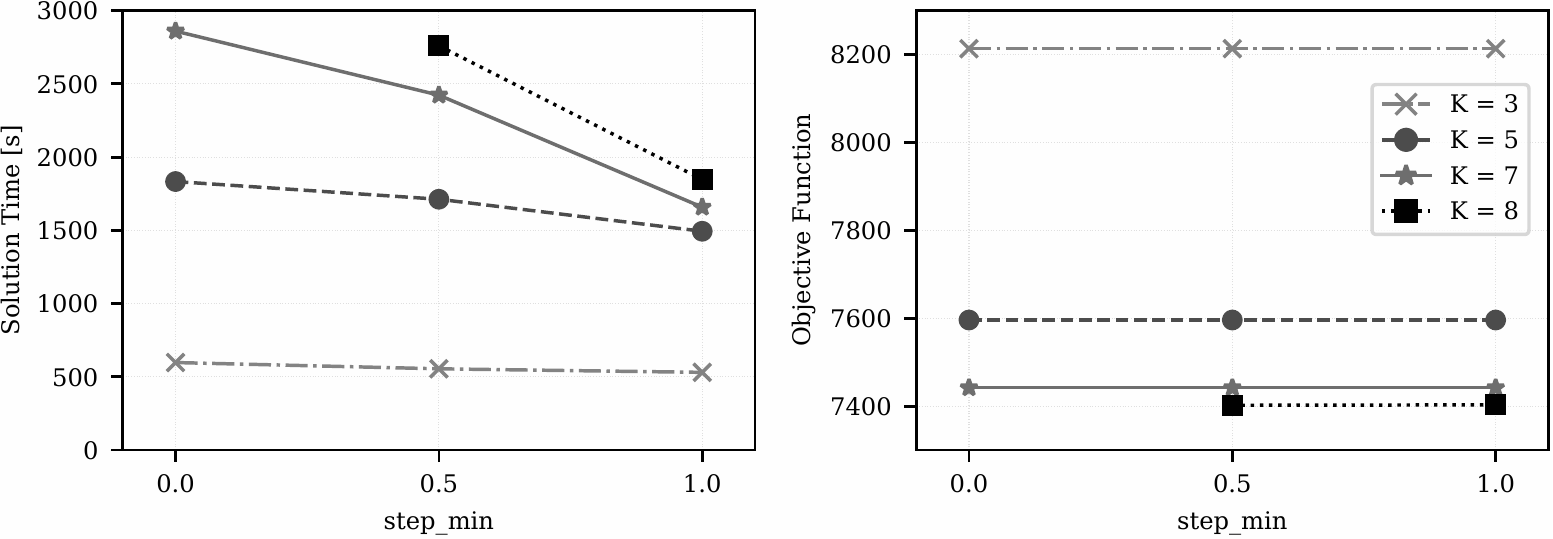}
      \caption{Objective function and solution times for different values of $K$ and \texttt{step\_min}.}
    
    \label{fig:experiment_step_min}
\end{figure}

On a different front, the solutions times T$^{\rm ISO}$ and T$^{\rm RLX}$ feature a steady increase as $K$ grows. This is consistent with the computational complexity of our algorithm, which depends linearly on $K$. Interestingly, T$^{\rm RLX}$ is substantially smaller than T$^{\rm ISO}$ for $K < 8$. The reason for this is that the lower bound we compute by relaxing the monotonicity constraint of the fitted curve naturally produces, however, a  fit that is non-increasing and thus, globally optimal. In contrast, when $K \geq 8$, such a lower bound does no longer coincide with the globally optimal solution and as a result,  T$^{\rm RLX}$ ends up surpassing the time limit set to one day (which is also exceeded by T$^{\rm ISO}$). To support this argument, Figure \ref{fig:realistic_k78} shows the fitted curves provided by the monotonicity-relaxed lower bound for $K=7$ and $K=8$. Notice that, if $K = 7$, the fitted curve associated with this lower bound is non-increasing, which allows our algorithm to certificate that this curve is, in fact, the optimal one in around 2900 seconds (with essentially all that time devoted to computing such a lower bound, logically). In contrast, when $K=8$, the curve delivered by the monotonicity-relaxed lower bound features a tiny step that destroys its otherwise non-increasing appearance. It is clear that this tiny step can only be attributed to the random nature of the data and not to the price-sensitivity of the pool of flexible consumers. Indeed, it is not reasonable to expect that a price variation lower than \euro{0.1}/MWh has such an impact on the consumption of the pool. In order to discard these implausible non-monotone stepwise functions from the family $\mathcal{F}$, our algorithm also includes the possibility to enforce a minimum arc length, that is, a minimum step size. Very conveniently, besides, this constraint helps reduce the solution time of our algorithm by pruning some paths in the graph that become thus infeasible and by increasing the chances that the monotonicity-relaxed lower bound corresponds to the optimal curve fitting. To illustrate the impact of the constraint on the minimum arc length on the computational performance of the proposed algorithm, we provide Figure~\ref{fig:experiment_step_min}, which shows the solution time for various step sizes and number $K$ of arcs. It can be seen that, while the case $K=8$ cannot be solved to optimality within a day time, if no constraint on the minimum arc length is enforced, the solution time is drastically reduced below 3000 seconds when a (very small) minimum block size of \euro0.5/MWh is imposed.

\section{Conclusions}\label{Conc}

In this paper, we have developed an algorithm to compute the curve that best fits to a certain cloud of data points in the sense of the least square error, under the conditions that said curve must be monotone and stepwise with a maximum number of steps. The proposed algorithm has been shown to run in polynomial time and is based on the finding that the curve-fitting problem can be addressed as a shortest-path type of problem. We have also proposed several strategies to cut down the execution time of the algorithm, all of which are based on computing upper and lower bounds that reduce the number of paths that the algorithm needs to explore. More specifically, the upper bound is given by a feasible solution that is swiftly built by combining the isotonic fit with clustering. The relaxation of either the constraint on the maximum number of steps or the monotonicity condition provides two different lower bounds, with the former being computationally much cheaper than the latter and not always necessarily looser. Our algorithm also allows for setting a minimum step length. This constraint notoriously speeds up the algorithm by pruning infeasible paths, while avoiding curve fits with implausible spurious tiny steps. 

Through a series of numerical experiments built on both synthetic and realistic data sets, we have demonstrated that our algorithm, in conjunction with the proposed acceleration strategies, is robust to the level of noise in the data and able to certificate the globally optimal curve in less than a few hours for sample sizes in the order of the thousands of data points. Furthermore, through a data set comprising  power-price  measurements  at  the  main  substation  of  a  distribution power  grid, we have shown that our algorithm serves as an useful tool to estimate the bidding curve whereby the distributed energy sources in the grid can trade in wholesale energy markets.
The extension of our algorithm to a multivariate setup is clearly an avenue of potentially fruitful research.

%
%
%


\section*{Acknowledgments}
This research has received funding from the European Research Council (ERC) under the European Union's Horizon 2020 research and innovation programme (grant agreement no. 755705). This work was also supported in part by the Spanish Ministry of Economy, Industry and Competitiveness and the European Regional Development Fund (ERDF) through project ENE2017-83775-P. Martine Labb\'e has been partially supported by the Fonds de la
Recherche Scientifique - FNRS under Grant(s) no PDR T0098.18.

\section*{CRediT author statement}

\textbf{V\'ictor Bucarey}: Conceptualization, Methodology, Formal Analysis, Validation, Writing - Original Draft, Software. \textbf{Martine Labb\'e}: Conceptualization, Methodology, Formal Analysis, Writing - Original Draft. \textbf{Juan Miguel Morales}: Conceptualization, Methodology, Formal Analysis, Validation, Writing - Original Draft, Supervision.  \textbf{Salvador Pineda}:  Methodology, Formal Analysis, Validation, Software, Data curation. 


\bibliographystyle{apalike} 
\bibliography{References} 



\end{document}